\documentclass{amsart}

\usepackage{amssymb}
\usepackage{color}

\newcommand{\mfT}{{\mathfrak T}}
\newcommand{\R}{\mathbb R}
\newcommand{\N}{\mathbb N}

\newcommand{\dsp}{\displaystyle}
\newtheorem{proposition}{Proposition}
\newtheorem{theorem}{Theorem}
\newtheorem{definition}{Definition}
\newtheorem{remark}{Remark}

\newtheorem{lemma}{Lemma}

\title[Large Time existence For 1D Green-Naghdi equations]{Large Time existence For 1D Green-Naghdi equations}
\author{Samer Israwi}
\address{Universit\'e Bordeaux I; IMB, 351 Cours de la Lib\'eration, 33405 Talence Cedex, France}
\thanks{Ce travail a b\'en\'efici\'e d'une aide de l'Agence Nationale de la Recherche portant  la r\'ef\'erence ANR-08-BLAN-0301-01}
\email{Samer.Israwi@math.u-bordeaux1.fr}

\begin{document}

\maketitle
\begin{abstract}
We consider here the $1D $ Green-Naghdi  equations that are commonly used in coastal oceanography  to
describe the propagation of large amplitude surface waves.
We show that  the solution of the  Green-Naghdi equations  can be constructed by a 
standard Picard iterative scheme so that there is no loss of regularity of the solution
with respect to the initial condition.  
\end{abstract}
\section{Introduction}
\subsection{Presentation of the problem}
The  water-waves problem for an ideal liquid
consists of describing the motion of the free surface and the evolution of the velocity
field of a layer of perfect, incompressible, irrotational fluid under the influence of
gravity. This motion is described by the free surface Euler equations that are known to be well-posed after the works of Nalimov \cite{nalimov}, Yasihara \cite{yoshihara}, 
Craig \cite{craig}, Wu \cite{wu1, wu2} and Lannes \cite{lannes}. But, because of the 
complexity of these equations,  they are often replaced for pratical purposes by approximate asymptotic systems. 
The most prominent examples are the Green-Naghdi equations (GN) -- which is 
a widely used model in coastal oceanography (\cite{GN,BMP,GLN}  and, for instance, 
\cite{JKSGR,KBEW})--, the Shallow-Water equations, and the
Boussinesq systems; their range of validity depends on the physical characteristics
of the flow under consideration. In other words, they depend on certain assumptions made 
on the dimensionless parameters $\varepsilon$, $\mu$ defined as:
\begin{equation*}
\varepsilon=\frac{a}{h_0},\quad \mu=\frac{h_0^2}{\lambda^2};
\end{equation*}
where $a$ is the order of amplitude of the waves and the bottom variations;
$\lambda$ is the wave-length of the waves and the bottom variations; 
$h_0$ is the reference  depth. The  parameter $\varepsilon$ is  often called  nonlinearity parameter;
while  $\mu$   is  the  shallowness  parameter. In the shallow-water scaling  $(\mu\ll1)$, and without smallness assumption on  $\varepsilon$  one  can  derive 
the so-called Green-Naghdi equations (see \cite{GN,LB} for a derivation and \cite{AL} for a
rigorous justification) also called Serre 
or fully nonlinear Boussinesq equations
 \cite{MB}.

 In nondimensionalized variables,  denoting by 
$\zeta(t,x)$ and $u(t,x)$ the parameterization of the surface  and the vertically averaged
horizontal component of the velocity at time $t$, and by $b(x)$ the parameterization
of the bottom, the equations read
 \begin{equation}
\left\{
\begin{array}{lc}
\dsp\partial_t\zeta+\nabla\cdot (hu)=0,\vspace{1mm}\\
\dsp(h+\mu h\mathcal{T}[h,\varepsilon b])\partial_t u+h\nabla\zeta+\varepsilon h( u\cdot \nabla )u\vspace{1mm}\\
\dsp\indent+\mu\varepsilon\big\lbrace-\frac{1}{3}\nabla[(h^3((u\cdot \nabla)(\nabla \cdot u)-(\nabla\cdot u)^2)]+ h\Re[h,\varepsilon b]u\big\rbrace=0,
\end{array}
\right.
\label{GN}
\end{equation}
where $h=1+\varepsilon(\zeta- b)$ and 
\begin{eqnarray*}
\mathcal{T}[h,\varepsilon b]W
=-\frac{1}{3h}\nabla(h^3\nabla\cdot W)+\frac{\varepsilon}{2h}[\nabla(h^2\nabla b\cdot W)-h^2 \nabla b\nabla\cdot W]+\varepsilon^2\nabla b \nabla b\cdot W,
\end{eqnarray*}
while the purely topographical term $\Re[h,\varepsilon b]u$ is defined as:
\begin{eqnarray*}
\Re[h,\varepsilon b]u&=&\frac{\varepsilon}{2h}[\nabla(h^2(u\cdot \nabla)^2b)-h^2((u\cdot \nabla)(\nabla \cdot u)-(\nabla\cdot u)^2)\nabla b]\\&&+\varepsilon^2((u\cdot \nabla)^2b)\nabla b.
\end{eqnarray*}
This model is often used in coastal oceanography because it takes into account the dispersive 
 effects neglected by the shallow-water and it is more nonlinear than the Boussinesq equations. 
A recent  rigorous justification of the GN model was given by Li \cite{li} in $1D$ and for flat bottoms, and by B. Alvarez-Samaniego and D. Lannes  \cite{AL}
in 2008 in the general case. This latter reference relies on  well-posedness results for these equations given in \cite{AL2} and based on general well-posedness results 
for evolution equations using a Nash-Moser scheme. The result of \cite{AL2} covers both the case of $1D$ and $2D$
surfaces, and allows for non flat bottoms. The reason why  a Nash-Moser scheme is used there is because
the estimates on the linearized equations exhibit losses of derivatives. However, in the $1D$ case with flat bottoms, 
such losses do not occur and it is possible to construct a solution with a standard Picard iterative scheme
as in \cite{li}. Our goal here is to show that it is also possible to use such a simple scheme in the $1D$ case with non flat bottoms, 
thanks to a careful analysis of the linearized equations. 
\subsection{Organization of the paper}
We start by giving some preliminary results  in Section \ref{pr};  the main theorem is
then stated in Section \ref{la} and proved in Section \ref{mr}. 
Finally, in Appendix \ref{appendix}, we give the existence and uniqueness of a 
solution to the linear Cauchy problem associated to the Green-Naghdi equations. The proof of the energy conservation, stated in the main theorem,  is given in Appendix \ref{sectconserv}.
\subsection{Notation}
We denote by $C(\lambda_1, \lambda_2,...)$ a constant depending on the parameters 
$\lambda_1$, $\lambda_2$, ... and \emph{whose dependence on the $\lambda_j$ is always assumed to be nondecreasing}.\\
The notation $a\lesssim b$ means that $a\leq Cb$, for some nonegative constant $C$ whose exact expression is of no importance (\emph{in particular, it is independent of the small parameters involved}).\\
Let $p$ be any constant
with $1\leq p< \infty$ and denote $L^p=L^p(\R)$ the space of all Lebesgue-measurable functions 
$f$ with the standard norm $$\vert f \vert_{L^p}=\big(\int_{\R}\vert f(x)\vert dx\big)^{1/p}<\infty.$$ When $p=2$,
we denote the norm $\vert\cdot\vert_{L^2}$ simply by $\vert\cdot\vert_2$. The inner product of any functions $f_1$
and $f_2$ in the Hilbert space $L^2(\R)$ is denoted by
$$
(f_1,f_2)=\int_{\R}f_1(x)f_2(x) dx.
$$  
The space $L^\infty=L^\infty(\R)$ consists of all essentially bounded, Lebesgue-measurable functions
$f$ with the norm
$$
\vert f\vert_{L^\infty}= \hbox{ess}\sup \vert f(x)\vert<\infty.
$$
We denote by $W^{1,\infty}=W^{1,\infty}(\R)=\{f\in L^\infty, \partial_x f\in L^{\infty}\}$ endowed with its canonical norm.\\
For any real constant $s$, $H^s=H^s(\R)$ denotes the Sobolev space of all tempered
distributions $f$ with the norm $\vert f\vert_{H^s}=\vert \Lambda^s f\vert_2 < \infty$, where $\Lambda$ 
is the pseudo-differential operator $\Lambda=(1-\partial_x^2)^{1/2}$.\\
For any functions $u=u(x,t)$ and $v(x,t)$
defined on $\R\times [0,T)$ with $T>0$, we denote the inner product, the $L^p$-norm and especially
the $L^2$-norm, as well as the Sobolev norm,
with respect to the spatial variable $x$, by $(u,v)=(u(\cdot,t),v(\cdot,t))$, $\vert u \vert_{L^p}=\vert u(\cdot,t)\vert_{L^p}$, 
$\vert u \vert_{L^2}=\vert u(\cdot,t)\vert_{L^2}$ , and $ \vert u \vert_{H^s}=\vert u(\cdot,t)\vert_{H^s}$, respectively.\\
Let $C^k(\R)$ denote the space of  
$k$-times continuously differentiable functions and $C^{\infty}_0(\R)$ denote the space of infinitely differentiable
functions, with compact support in $\R$; we also denote by $C^\infty_b(\R)$ the space of infinitely differentiable functions that are bounded together with all their derivatives.\\
Let $f$ be a function of the independent variables
$x_1$, $x_2$,...,$x_m$; its partial derivative with respect to $x_k$ is denoted by 
$\partial_{x_k}f=f_{x_k}$ for $1 \leq k \leq m$.\\
For any closed operator $T$ defined on a Banach space $X$ of functions, the commutator $[T,f]$ is defined
 by $[T,f]g=T(fg)-fT(g)$ with $f$, $g$ and $fg$ belonging to the domain of $T$.\\
\section{Well-posedness of the Green-Naghdi equations in $1D$ }\label{wp}

For one dimensional surfaces, the Green-Naghdi equations (\ref{GN}) can be simplified,
after some computations, 
into 
\begin{equation}
\left\{
\begin{array}{lc}
\dsp\partial_t\zeta+\partial_x(hu)=0,\vspace{1mm}\\
\dsp(h+\mu h\mathcal{T}[h,\varepsilon b])[\partial_t u+\varepsilon u\partial_xu]+h\partial_x\zeta
+\varepsilon\mu hQ[h,\varepsilon b](u)=0\vspace{1mm}\\
\end{array}
\right.
\label{GN1}
\end{equation}
where $h=1+\varepsilon(\zeta- b)$ and
$$
\mathcal{T}[h,\varepsilon b]w
=-\frac{1}{3h}\partial_x(h^3 w_x)+\frac{\varepsilon}{2h}[\partial_x(h^2 b_x w)-h^2 b_x w_x]+\varepsilon^2b_x^2  w,
$$
$$
Q[h,\varepsilon b](w)=\frac{2}{3h}\partial_x (h^3 w_x^2)+\varepsilon h w_x^2 b_x+\varepsilon\frac{1}{2h}\partial_x (h^2 w^2 b_{xx})+\varepsilon^2 w^2 b_{xx}b_x. 
$$
\begin{remark}
The interest of the formulation (\ref{GN1}) of the Green-Naghdi equation is that all
the third order derivatives of $u$ have been factorized by $(h+\mu h\mathcal{T}[h,\varepsilon b])$. Indeed, $Q[h,\varepsilon b]$ is a second order differential operator. This was used in \cite{li} in the case of flat bottoms ($b=0$).
\end{remark}
\subsection{Preliminary results}\label{pr}
For the sake of simplicity, we write $$\mfT=h+\mu h\mathcal{T}[h,\varepsilon b].$$ We always assume that  the 
nonzero depth condition 
\begin{equation}\label{depthcond}
\exists\; h_0 >0, \quad \inf_{x\in \R} h\ge h_0,\quad h=1+\varepsilon(\zeta- b)
\end{equation}
is valid  initially, which is a necessary condition for the GN system (\ref{GN1}) to be  physically valid.
We shall demonstrate that the operator $\mfT$ plays an important role in the energy estimate and the local well-posedness of the GN system 
(\ref{GN1}). Therefore, we give here some of its properties.

The following lemma gives an important invertibility result on $\mfT$.
\begin{lemma}\label{proprim}
Let $b\in C_b^{\infty}(\R)$ and $\zeta \in W^{1,\infty}(\R)$ be such that (\ref{depthcond}) is satisfied.
 Then the operator 
 $$
\mfT: H^2(\R)\longrightarrow L^2(\R)
$$ 
is well defined, one-to-one and onto.
\end{lemma}
\begin{remark}\label{remark b}
Here and throughout the rest of this paper, and for the sake of simplicity, we do not try to give some optimal regularity assumption on 
the bottom parameterization $b$. This could easily be done, but is of no interest for
our present purpose. Consequently, we ommit to write the dependance on $b$ of the
different quantities that appear in the proof. 
\end{remark}
\begin{proof}
In order to prove the invertibility of $\mfT$, let us
first remark that  the quantity $\vert v\vert^2_{*}$ defined as
$$
\vert v\vert^2_{*}=\vert v\vert^2_2+\mu\vert \partial_x v\vert^2_{2},
$$
is equivalent to the $H^1(\R)$-norm but not uniformly with respect to $\mu\in (0,1)$.
 We define by $H^1_*(\R)$ the space $H^1(\R)$ endowed with this norm.
The bilinear form:
\begin{equation*}
a(u,v)=(hu,v)+\mu\big(h\big(\frac{h}{\sqrt{3}}u_x-\frac{\sqrt{3}}{2}\varepsilon b_xu\big),\frac{h}{\sqrt{3}}v_x-\frac{\sqrt{3}}{2}\varepsilon b_xv\big)+\frac{\mu\varepsilon^2}{4}(hb_xu,b_xv).
\end{equation*}
is obviously continous on $H^1_*(\R)\times H^1_*(\R)$. Remarking that 
\begin{eqnarray*}
a(v,v)=(\mfT v,v)&=&(hv,v) \\
&&+\mu\big(h\big(\frac{h}{\sqrt{3}}v_x-\frac{\sqrt{3}}{2}\varepsilon b_xv\big),\frac{h}{\sqrt{3}}v_x-\frac{\sqrt{3}}{2}\varepsilon b_xv\big)+\frac{\mu\varepsilon^2}{4}(hb_xv,b_xv),
\end{eqnarray*}
 we have
\begin{eqnarray*}
 \vert v\vert^2_{*}&\leq&\vert v\vert^2_2+\frac{3\mu}{h_0^2}\vert \frac{h}{\sqrt{3}} v_x\vert^2_2\\
&\leq& \vert v\vert^2_2+\frac{6\mu}{h_0^2}\Big(\vert \frac{h}{\sqrt{3}} v_x-\frac{\sqrt{3}}{2}\varepsilon
b_xv\vert^2_2+\frac{3\varepsilon^2}{4} \vert b_xv\vert^2_2\Big).
\end{eqnarray*}
One deduces that
$$\max\Big\{1,\frac{18}{h_0^2}\Big\}\Big(\vert v \vert^2_2
+\mu\vert \frac{h}{\sqrt{3}}v_x-\frac{\sqrt{3}}{2}\varepsilon b_xv \vert^2_2+ \frac{\mu\varepsilon^2}{4} \vert b_xv\vert^2_2\Big)\ge \vert v\vert^2_{*}.$$
Since from  (\ref{depthcond}) we also get
\begin{eqnarray*}
a(v,v)&\ge& h_0\vert v \vert^2_2
+\mu h_0\Big(\vert \frac{h}{\sqrt{3}}v_x-\frac{\sqrt{3}}{2}\varepsilon b_xv \vert^2_2+ \frac{\mu\varepsilon^2}{4} \vert b_xv\vert^2_2\Big),
\end{eqnarray*}
it is easy to deduce that 
\begin{eqnarray}\label{minorantim}
a(v,v)&\ge& \frac{h_0}{\max\big\{1,\frac{18}{h_0^2}\big\}}\vert v\vert^2_{*}.
\end{eqnarray}
In particular, $a$ is coercive on $H^1_{*}$.
Using Lax-Milgram lemma, for all $f \in L^2(\R)$,
there exists unique $u\in H^1_*(\R)$
such that, for all $v\in H^1_{*}(\R)$
\begin{equation*}
a(u,v)=(f,v);
\end{equation*}
equivalently, there is a unique variational solution to the equation 
\begin{equation}\label{reeq}
\mfT u=f.
\end{equation}
We then get from the definition of $\mfT$ that
$$
\partial_x^2u=\frac{hu+\frac{\varepsilon\mu}{2}\partial_x(h^2b_x)u+\varepsilon^2\mu hb^2_x u-\frac{\mu}{3}\partial_xh^3\partial_xu-f}{\frac{\mu h^3}{3}};
$$
since $u\in H^1(\R)$ and $f\in L^2(\R)$, we get $\partial_x^2u \in L^2(\R)$ and thus $u\in H^2(\R)$.
\end{proof}
The following lemma then gives some properties of the inverse operator $\mfT^{-1}$.
\begin{lemma}\label{proprim'}
Let $b\in C_b^{\infty}(\R)$, $t_0> 1/2$  and $\zeta \in H^{t_0+1}(\R)$ 
be such that (\ref{depthcond}) is satisfied. Then:\\\\
\quad (i) $\forall 0\leq s\leq t_0+1$, $\quad\vert \mfT^{-1}f\vert_{H^s}+\sqrt{\mu}\vert \partial_x\mfT^{-1}f\vert_{H^s}\leq 
C(\frac{1}{h_0},\vert h-1 \vert_{H^{t_0+1}})\vert f\vert_{H^s}$;\\\\
\quad (ii)  $\forall 0\leq s\leq t_0+1$, $\quad\sqrt{\mu}\vert \mfT^{-1}\partial_x g\vert_{H^s}\leq 
C(\frac{1}{h_0},\vert h-1 \vert_{H^{t_0+1}})
\vert g\vert_{H^s}$;\\\\
\quad (iii) If  $s\geq t_0+1$ and $\zeta\in H^s(\R)$ then:\\
$$
\parallel \mfT^{-1}\parallel_{H^s(\R)\rightarrow H^s(\R)}+
\sqrt{\mu}\parallel\mfT^{-1}\partial_x\parallel_{H^s(\R)\rightarrow H^s(\R)}\leq c_s,
$$
where $c_s$ is a constant depending on $\frac{1}{h_0}$,  $\vert h-1 \vert_{H^{s}}$
and independent of ($\mu$,$\varepsilon$) $\in (0,1)^2$.
\end{lemma}
\begin{proof}
\textbf{Step 1}. We prove that if $u\in H^1_*(\R)$ solves 
$$
\mfT u=f+\sqrt{\mu}\partial_xg
$$
for $f,g \in L^2(\R)$, then one has 
$$
\vert u\vert_{H^1_*}\leq C\big(\frac{1}{h_0}\big)\big(\vert f\vert_2+\vert g\vert_2\big).
$$
Indeed, multiplying the equation by $u$ and integrating by parts, one gets, with the notations used 
in the proof of lemma \ref{proprim}
$$
a(u,u)\leq (f,u)-(g,\sqrt{\mu}\partial_x u).
$$
We thus get from the proof of Lemma \ref{proprim} and Cauchy-Schwarz inequality that
$$
\frac{h_0}{\max\{1,\frac{18}{h_0^2}\}}\vert u\vert^2_{H^1_*}\leq \vert f\vert_2\vert u\vert_2+\vert g\vert_2\vert u\vert_{H^1_*},
$$
and the result follows easily.\\
\textbf{Step 2}. We prove here that 
$\vert \mfT^{-1}f\vert_{H^s}+\sqrt{\mu}\vert \partial_x\mfT^{-1}f\vert_{H^s}\leq 
C(\frac{1}{h_0},\vert h-1 \vert_{H^{t_0+1}})\vert f\vert_{H^s}$.
Indeed, if $f\in H^s$ and $u=\mfT^{-1}f$ then $\mfT u=f$.
Applying $\Lambda^s$ to this identity, we get 
\begin{eqnarray*}
 \mfT(\Lambda^s u)&=& \Lambda^sf+[\mfT,\Lambda^s]u\\
&=&\tilde{f}+\sqrt{\mu}\partial_x\tilde{g},
\end{eqnarray*}
with,
$$
\tilde{f}=\Lambda^sf-[\Lambda^s,h]u+\frac{\varepsilon\mu}{2}[\Lambda^s,h^2b_x]u_x-\varepsilon^2\mu[\Lambda^s,h^2b_x]u,
$$
and
$$
\tilde{g}=\frac{\sqrt{\mu}}{3}[\Lambda^s,h^3]u_x-\frac{\varepsilon\sqrt{\mu}}{2}[\Lambda^s,h^2b_x]u.
$$
Now, one can deduce from the commutator estimate (see e.g Lemma 4.6 of \cite{AL2})
\begin{equation}\label{com1}
\vert[\Lambda^s,F]G\vert_2\lesssim \vert \nabla F\vert_{H^{t_0}}\vert G\vert_{H^{s-1}}
\end{equation}
that 
\begin{eqnarray*}
 \vert \tilde{f}\vert_2+ \vert \tilde{g}\vert_2
&\leq& \vert f\vert_{H^s}+C\big(\frac{1}{h_0},\vert h-1 \vert_{H^{t_0+1}}\big)\big(\vert u\vert_{H^{s-1}}+\sqrt{\mu}\vert \partial_xu\vert_{H^{s-1}}\big).
\end{eqnarray*}
One can use Step 1 and a continuous induction on $s$ to show that the inequality (i) holds for $0\leq s\leq t_0+1$.\\
\textbf{Step 3}. We prove here that 
$\sqrt{\mu}\vert \mfT^{-1}\partial_x g\vert_{H^s}\leq 
C(\frac{1}{h_0},\vert h-1 \vert_{H^{t_0+1}})
\vert g\vert_{H^s}$.
Indeed, if $g\in H^s$ and $u=\sqrt{\mu}\mfT^{-1}\partial_xg$ then $\mfT u=\sqrt{\mu}\partial_xg$
and thus $$\mfT(\Lambda^su)=\tilde{f}+\sqrt{\mu}\partial_x\tilde{g},$$
with,
$$
\tilde{f}=-[\Lambda^s,h]u+\frac{\varepsilon\mu}{2}[\Lambda^s,h^2b_x]u_x-\varepsilon^2\mu[\Lambda^s,h^2b_x]u,
$$
and
$$
\tilde{g}=\Lambda^sg+\frac{\sqrt{\mu}}{3}[\Lambda^s,h^3]u_x-\frac{\varepsilon\sqrt{\mu}}{2}[\Lambda^s,h^2b_x]u.
$$
Proceeding now as for the Step 2, one can deduce (ii).\\
\textbf{Step 4}. If $s\ge t_0+1$ then one can prove (iii) proceeding as in Step 2 and 3 above,
but replacing the commutator estimate (\ref{com1}) by the following one
\begin{equation}\label{com2}
\vert[\Lambda^s,F]G\vert_2\lesssim \vert \nabla F\vert_{H^{s-1}}\vert G\vert_{H^{s-1}}.
\end{equation}
\end{proof}
\subsection{Linear analysis}\label{la}
In order to rewrite the GN equations (\ref{GN1}) in a condensed form, let us decompose $Q[h,\varepsilon b](u)$ as
\begin{equation*}
\varepsilon\mu hQ[h,\varepsilon b](u)=Q_1[U]u_x+q_2(U)
\end{equation*}
where $U=(\zeta,u)^T$ and
\begin{eqnarray}\label{decompose}
Q_1[U]f&=&\frac{2}{3}\varepsilon\mu\partial_x(h^3u_xf)+\varepsilon^2\mu h^2b_xu_xf+\varepsilon^2\mu h^2b_{xx}uf \nonumber\\\\
q_2(U)&=&\varepsilon^3\mu hb_{xx}b_xu^2+\frac{1}{2}\varepsilon^2\mu \partial_x(h^2b_{xx})u^2.\nonumber
\end{eqnarray}
The Green-Naghdi equations (\ref{GN1}) can be written after applying $\mfT^{-1}$ to both sides of the second equation
in  (\ref{GN1}) as
\begin{equation}\label{condensedeq}
\partial_tU+A[U]\partial_xU+B(U)=0,
\end{equation}
with $U=(\zeta,u)^T$ and where\\

\begin{equation}
A[U]=\left(
\begin{array}{cc}
\varepsilon u&h\\\\
\mfT^{-1}(h\cdot)& \varepsilon u+\mfT^{-1}Q_1[U]
\end{array}
\right)
\end{equation}
and\\
\begin{equation}
B(U)=\left(
\begin{array}{c}
\varepsilon b_xu\\\\
 \mfT^{-1}q_2(U)
\end{array}
\right).
\end{equation}

This subsection is devoted to the proof of energy estimates for the following initial value problem around some reference state 
$\underline{U}=(\underline{\zeta},\underline{u})^T$:
\begin{equation}\label{GNlsys}
	\left\lbrace
	\begin{array}{l}
	\dsp\partial_t U+A[\underline{U}]\partial_x U+B (\underline{U})=0;
        \\
	\dsp U_{\vert_{t=0}}=U_0.
	\end{array}\right.
\end{equation}

We define now the $X^s$ spaces, which are the energy spaces for this problem.
\begin{definition}\label{defispace}
 For all $s\ge 0$ and $T>0$, we denote by $X^s$ the vector space $H^s(\R)\times H^{s+1}(\R)$ endowed with the norm
$$
\forall\; U=(\zeta,u) \in X^s, \quad \vert U\vert^2_{X^s}:=\vert \zeta\vert^2 _{H^s}+\vert u\vert^2 _{H^s}+ \mu\vert \partial_xu\vert^2 _{H^s},
$$
while $X^s_T$ stands for $C([0,\frac{T}{\varepsilon}];X^{s})$ endowed with its canonical norm.
\end{definition}
First remark that a symmetrizer for $A[\underline{U}]$ is given by
\begin{equation}
S=\left( 
\begin{array}{cc}
 1& 0 \\\\ 
0& \underline{\mfT}
\end{array}
\right),
\end{equation}\\
with 
$\underline{h}=1+\varepsilon(\underline{\zeta}-b)$ and
$\underline{\mfT}=\underline{h}+\mu\underline{h}\mathcal{T}[\underline{h},\varepsilon b].$
A natural energy for the IVP (\ref{GNlsys}) is given by
\begin{equation}\label{es}
 E^s(U)^2=(\Lambda^sU,S\Lambda^sU).
\end{equation}
The link between $ E^s(U)$ and the $X^s$-norm is investigated in the following Lemma.
\begin{lemma}\label{lemmaes}
Let $b\in C_b^{\infty}(\R)$,  $s\geq 0$ and $ \underline{\zeta}\in W^{1,\infty}(\R)$. Under the condition (\ref{depthcond}),
$E^s(U)$ is uniformly equivalent to the $\vert \cdot\vert_{X^s}$-norm  with respect to $(\mu, \varepsilon) \in (0,1)^2$:
$$
E^s(U) \leq C\big(\vert \underline{h}\vert_{L^{\infty}},\vert \underline{h}_x\vert_{L^{\infty}}\big)\vert U\vert_{X^s},
$$
and
$$
\vert U \vert_{X^s}\leq C\big(\frac{1}{h_0}\big) E^s(U).
$$
\end{lemma}
\begin{proof}
 Notice first that
$$
E^s(U)^2=\vert \Lambda^s\zeta\vert_{2}^2 +(\Lambda^su,\underline{\mfT}\Lambda^su),
$$
 one gets the first estimate using the explicit expression of $\underline{\mfT}$, integration by parts and Cauchy-Schwarz inequality. \\
The other inequality can be proved by using that $\inf_{x\in\R}h\ge h_0>0$ and proceeding as in the proof of  Lemma \ref{proprim}. 
\end{proof}
We prove now the energy estimates in the following proposition:
\begin{proposition}\label{ESprop}
Let  $b\in C_b^{\infty}(\R)$,  $t_0>1/2$, $s\geq t_0+1$. Let also $\underline{U}=(\underline{\zeta}, \underline{u})^T$ $\in X^{s}_{T}
$ be such that $\partial_t \underline{U} \in X^{s-1}_{T}$ 
and satisfying the condition  (\ref{depthcond}) on $[0,\frac{T}{\varepsilon}]$. Then for all   $U_0\in X^{s}$
there exists a unique solution  $U=(\zeta, u)^T$ $\in X^{s}_{T} $ to (\ref{GNlsys}) and for all $0\leq t\leq\frac{T}{\varepsilon}$
$$
E^s(U(t))\leq e^{\varepsilon\lambda_{T} t}E^s(U_0)+\varepsilon \int^{t}_{0} e^{\varepsilon\lambda_T( t-t')}C(E^s(\underline{U})(t'))dt'.
$$
For some $\lambda_{T}=\lambda_{T}(\sup_{0\leq t\leq T/\varepsilon}E^s(\underline{U}(t)),\sup_{0\leq t\leq T/\varepsilon}\vert\partial_t\underline{h}(t) \vert_{L^{\infty}})$ .
\end{proposition}
\begin{proof}
Existence and uniqueness of a solution to the IVP (\ref{GNlsys})  is achieved in ap- 
pendix A and we thus focus our attention on the proof 
of the energy estimate. For any $\lambda \in \R$, we compute
$$
e^{\varepsilon\lambda t}\partial_t(e^{-\varepsilon\lambda t}E^s(U)^2)=-\varepsilon\lambda E^s(U)^2 +\partial_t(E^s(U)^2).
$$
Since
$$
 E^s(U)^2=(\Lambda^sU,S\Lambda^sU),
$$
we have
\begin{equation}
\partial_t(E^s(U)^2)=2(\Lambda^s\zeta,\Lambda^s\zeta_t)+2(\Lambda^su,\underline{\mfT}\Lambda^su_t)+(\Lambda^su,[\partial_t,\underline{\mfT}]\Lambda^su).
\end{equation}
One gets  using the equations  (\ref{GNlsys}) and integrating by parts,
\begin{eqnarray}\label{qtsctr}
\frac{1}{2}e^{\varepsilon\lambda t}\partial_t(e^{-\varepsilon\lambda t}E^s(U)^2)&=&-\frac{\varepsilon\lambda}{2}E^s(U)^2 
-(SA[\underline{U}]\Lambda^s\partial_x U,\Lambda^s U) \nonumber\\\nonumber\\&&- \big(\big[\Lambda^s,A[\underline{U}]\big]\partial_x U,S\Lambda^s U\big)
-(\Lambda^s B(\underline{U}),S\Lambda^s U)\nonumber\\\nonumber\\&&+\frac{1}{2}(\Lambda^su,[\partial_t,\underline{\mfT}]\Lambda^su).
\end{eqnarray}
We now turn to bound from above the different components of the r.h.s of (\ref{qtsctr}).\\
$\bullet$ Estimate of $(SA[\underline{U}]\Lambda^s\partial_x U,\Lambda^s U)$.
 Remarking  that
 \begin{equation*}
SA[\underline{U}]=\left(
\begin{array}{cc}
\varepsilon \underline{u}&\underline{h}\\\\
\underline{h}& \underline{\mfT}(\varepsilon \underline{u}\cdot)+Q_1[\underline{U}]
\end{array}
\right),
\end{equation*}
we get
\begin{eqnarray*}
(SA[\underline{U}]\Lambda^s\partial_x U,\Lambda^s U)&=&  (\varepsilon\underline{u}\Lambda^s\zeta_x,\Lambda^s\zeta)+
(\underline{h}\Lambda^su_x,\Lambda^s\zeta)\\\\&&
+(\underline{h}\Lambda^s \zeta_x,\Lambda^su)+\big((\underline{\mfT}(\varepsilon\underline{u}\cdot)+Q_1[\underline{U}])\Lambda^s u_x,\Lambda^su\big)\\\\&&
=: A_1+A_2+A_3+A_4.
\end{eqnarray*}
We now focus to control $(A_j)_{1\leq j\leq 4}$.\\
$-$ Control of $A_1$. Integrating by parts, one obtains 
$$
A_1= (\varepsilon\underline{u}\Lambda^s\zeta,\Lambda^s\zeta_x)=-\frac{1}{2}(\varepsilon\underline{u}_x\Lambda^s\zeta,\Lambda^s\zeta)
$$
one can conclude by Cauchy-Schwarz inequality that
$$
\vert A_1\vert \leq\varepsilon C(\vert \underline{u}_x\vert_{L^{\infty}})E^s(U)^2.
$$
$-$ Control of $A_2+A_3$. First remark that 
$$
\vert A_2+A_3\vert=\vert(\underline{h}_x\Lambda^s u,\Lambda^s\zeta)\vert\leq \vert\underline{h}_x \vert_{L^{\infty}}E^s(U)^2;
$$
we get,
$$
\vert A_2+A_3\vert\leq\varepsilon C(\vert\underline{h}_x\vert_{L^\infty})E^s(U)^2.
$$
$-$ Control of $A_4$. One computes,
\begin{eqnarray*}
A_4&=&\varepsilon\big(\underline{\mfT} (\underline{u}\Lambda^s u_x),\Lambda^su\big)+(Q_1[\underline{U}]\Lambda^s u_x,\Lambda^su)\\&&
=:A_{41}+A_{42}.
\end{eqnarray*}
 Note that
\begin{eqnarray*}
A_{41}&=&\varepsilon(\underline{h} \; \underline{u}\Lambda^s u_x,\Lambda^su)+\frac{\varepsilon\mu}{3}(\underline{h}^3 \; (\underline{u}\Lambda^s u_x)_x,\Lambda^su_x)\\\\&&
-\frac{\varepsilon^2\mu}{2}(\underline{h}^2 \; b_x(\underline{u} \Lambda^s u_x)_x,\Lambda^su)-\frac{\varepsilon^2\mu}{2}(\underline{h}^2 b_x \underline{u}\Lambda^s u_x,\Lambda^su_x)\\\\&&
+\varepsilon^3\mu(\underline{h} \; \underline{u} b_x^2\Lambda^s u_x,\Lambda^su);
\end{eqnarray*}
since
\begin{eqnarray*}
(\underline{h}^3 \; (\underline{u}\Lambda^s u_x)_x,\Lambda^su_x)&=& \dsp\frac{1}{2}\big(-(\underline{h}^3_x \; \underline{u}\Lambda^s u_x,\Lambda^su_x)+(\underline{h}^3 \; \underline{u}_x\Lambda^s u_x,\Lambda^su_x)\big),
\end{eqnarray*}
by using successively integration by parts  and the Cauchy-Schwarz inequality, one obtains directly:
$$
\vert A_{41}\vert \leq \varepsilon C(\vert \underline{u}\vert_{W^{1,\infty}},\vert \underline{\zeta}\vert_{W^{1,\infty}})E^s(U)^2.
$$
For $A_{42}$, remark that 
\begin{eqnarray*}
\vert A_{42}\vert&=&\vert(Q_1[\underline{U}]\Lambda^s u_x,\Lambda^su) \vert\\\\
&=& \big\vert -\frac{2}{3}\varepsilon\mu(\underline{h}^3\underline{u}_x\Lambda^s u_x,\Lambda^su_x)+\varepsilon^2\mu(\underline{h}^2 \; \underline{u}_x b_x\Lambda^s u_x,\Lambda^su)\\\\&&
\quad+\varepsilon^2\mu(\underline{h}^2 \; \underline{u} b_{xx}\Lambda^s u_x,\Lambda^su)\big\vert
\end{eqnarray*}
therefore
$$
\vert A_{42}\vert \leq \varepsilon C(\vert \underline{u}\vert_{W^{1,\infty}},\vert \underline{\zeta}\vert_{W^{1,\infty}})E^s(U)^2.
$$
This shows that 
$$
\vert A_{4}\vert \leq \varepsilon C(\vert \underline{u}\vert_{W^{1,\infty}},\vert \underline{\zeta}\vert_{W^{1,\infty}})E^s(U)^2.
$$
$\bullet$ Estimate of $ \big(\big[\Lambda^s,A[\underline{U}]\big]\partial_x U,S\Lambda^s U\big)$.
 Remark first that
 \begin{eqnarray*}
\big(\big[\Lambda^s,A[\underline{U}]\big]\partial_x U,S\Lambda^s U\big)&=&([\Lambda^s,\varepsilon \underline{u}]\zeta_x,\Lambda^s\zeta)+([\Lambda^s, \underline{h}]u_x,\Lambda^s\zeta)\\\\
&&+([\Lambda^s, \underline{\mfT}^{-1}\;\underline{h}]\zeta_x,\underline{\mfT}\Lambda^su)+([\Lambda^s, \varepsilon \underline{u}]u_x,\underline{\mfT}\Lambda^su)\\\\&&
+\big(\big[\Lambda^s, \underline{\mfT}^{-1}Q_1[\underline{U}]\big]u_x,\underline{\mfT}\Lambda^su\big)\\\\
&=:& B_1+B_2+B_3+B_4+B_5.
\end{eqnarray*}
$-$ Control of $B_1+B_2=([\Lambda^s,\varepsilon \underline{u}]\zeta_x,\Lambda^s\zeta)+([\Lambda^s, \underline{h}]u_x,\Lambda^s\zeta)$.
Since $s\geq t_0+1$, we can use the commutator estimate (\ref{com2}) to get 
$$
\vert B_1+B_2\vert  \leq \varepsilon C(E^s(\underline{U}))E^s(U)^2.
$$\\ 
$-$ Control of $B_4=([\Lambda^s, \varepsilon \underline{u}]u_x,\underline{\mfT}\Lambda^su)$. By using  the explicit expression of  $\underline{\mfT}$ we get
\begin{eqnarray*}
B_4&=&([\Lambda^s, \varepsilon \underline{u}]u_x,\underline{h}\Lambda^s u)+\frac{\mu}{3}(\partial_x[\Lambda^s, \varepsilon \underline{u}]u_x,\underline{h}^3 \; \Lambda^s u_x)\\\\&&
-\frac{\varepsilon\mu}{2}([\Lambda^s, \varepsilon \underline{u}]u_x,\underline{h}^2 \;  b_x\Lambda^s u_x)
+\frac{\varepsilon\mu}{2}([\Lambda^s, \varepsilon \underline{u}]u_x,\partial_x(\underline{h}^2 \; b_x\Lambda^s u))
\\\\&&
+\varepsilon^2\mu([\Lambda^s, \varepsilon \underline{u}]u_x,\underline{h} \; b_x^2\Lambda^s u),
\end{eqnarray*}
using  the Cauchy-Schwarz inequality and the fact that 
$$
\partial_x[\Lambda^s, f]g=[\Lambda^s, f_x]g+[\Lambda^s, f]g_x
$$
one obtains directly:
$$
\vert B_4\vert  \leq \varepsilon C(E^s(\underline{U}))E^s(U)^2.
$$\\ 
$-$ Control of $B_3=([\Lambda^s, \underline{\mfT}^{-1}\;\underline{h}]\zeta_x,\underline{\mfT}\Lambda^su)$. Remark first that
\begin{equation*}
\underline{\mfT}[\Lambda^s, \underline{\mfT}^{-1}]\underline{h}\zeta_x=\underline{\mfT}[\Lambda^s, \underline{\mfT}^{-1}\;\underline{h}]\zeta_x
-[\Lambda^s,\underline{h}]\zeta_x;
\end{equation*}
morever, since $[\Lambda^s,\underline{\mfT}^{-1}]=-\underline{\mfT}^{-1}[\Lambda^s,\underline{\mfT}]\underline{\mfT}^{-1}$, one gets 
\begin{equation*}
\underline{\mfT}[\Lambda^s, \underline{\mfT}^{-1}\;\underline{h}]\zeta_x=-[\Lambda^s, \underline{\mfT}] \underline{\mfT}^{-1}\underline{h}\zeta_x
+[\Lambda^s,\underline{h}]\zeta_x,
\end{equation*}
and one can check by using  the explicit expression of  $\underline{\mfT}$ that
\begin{eqnarray*}
\underline{\mfT}[\Lambda^s, \underline{\mfT}^{-1}\;\underline{h}]\zeta_x&=& -[\Lambda^s, \underline{h}] \underline{\mfT}^{-1}\underline{h}\zeta_x
+\frac{\mu}{3}\partial_x\{[\Lambda^s, \underline{h}^3]\partial_x(\underline{\mfT}^{-1}\underline{h}\zeta_x)\} \\\\&&
-\frac{\varepsilon\mu}{2}
\partial_x[\Lambda^s, \underline{h}^2b_x] \underline{\mfT}^{-1}\underline{h}\zeta_x
+\frac{\varepsilon\mu}{2}
[\Lambda^s, \underline{h}^2b_x] \partial_x\underline{\mfT}^{-1}\underline{h}\zeta_x \\\\ &&
-\varepsilon^2\mu[\Lambda^s, \underline{h}b_x^2] \underline{\mfT}^{-1}\underline{h}\zeta_x
+[\Lambda^s,\underline{h}]\zeta_x.
\end{eqnarray*}
One deduces directly from Lemma \ref{proprim'}, an integration by parts, and Cauchy-Schwarz inequality 
that
\begin{eqnarray*}
\vert B_3\vert & \leq& C\big(\frac{1}{h_0},\vert \underline{h}-1 \vert_{H^{s}}\big)\; \Big\{\Big( \vert\underline{h}_x \vert_{H^{s-1}}+\frac{\varepsilon\mu}{2}\vert\underline{h}^2b_x \vert_{H^s}
+\varepsilon^2\mu\vert\underline{h}b_x^2 \vert_{H^s} \Big)\vert\underline{h}  \zeta_x\vert_{H^{s-1}}\\\\&&
+ \Big(\frac{\sqrt{\mu}}{3}\vert\underline{h}^3_x\vert_{H^{s-1}}+ \frac{\varepsilon\sqrt{\mu}}{2}\vert\underline{h}^2b_x \vert_{H^s}\Big)
\vert\underline{h} \zeta_x\vert_{H^{s-1}} 
 +\vert\underline{h}_x \vert_{H^{s-1}} \vert \zeta_x\vert_{H^{s-1}}\Big)\Big\} \vert \Lambda^su\vert_{H^1_*}.
\end{eqnarray*}
Finally,  since
$$
  \vert\underline{h}  \zeta_x\vert_{H^{s-1}}\leq C(E^s(\underline{U}))E^s(U)\quad\hbox{and}\quad
 \vert\underline{h} b_x^2\vert_{H^{s}}+\vert\underline{h}^2 b_x\vert_{H^{s}}\leq C( E^s(\underline{U})),
$$
we deduce
$$
\vert B_3\vert  \leq \varepsilon C(E^s(\underline{U}))E^s(U)^2.
$$\\ 
$-$ Control of $B_5=\big(\big[\Lambda^s, \underline{\mfT}^{-1}Q_1[\underline{U}]\big]u_x,\underline{\mfT}\Lambda^su\big)$. Let us first write 
\begin{equation*}
\underline{\mfT}\big[\Lambda^s, \underline{\mfT}^{-1}Q_1[\underline{U}]\big]u_x=-[\Lambda^s, \underline{\mfT}] \underline{\mfT}^{-1}Q_1[\underline{U}]u_x
+\big[\Lambda^s,Q_1[\underline{U}]\big]u_x
\end{equation*}
so, that 
\begin{eqnarray*}
\underline{\mfT}\big[\Lambda^s, \underline{\mfT}^{-1}Q_1[\underline{U}]\big]u_x&=& -[\Lambda^s, \underline{h}] \underline{\mfT}^{-1}Q_1[\underline{U}]u_x
+\frac{\mu}{3}\partial_x\{[\Lambda^s, \underline{h}^3]\partial_x(\underline{\mfT}^{-1}Q_1[\underline{U}]u_x)\} \\\\&&
- \frac{\varepsilon\mu}{2}
\partial_x\{[\Lambda^s, \underline{h}^2b_x] \underline{\mfT}^{-1}Q_1[\underline{U}]u_x\}+
\frac{\varepsilon\mu}{2}
[\Lambda^s, \underline{h}^2b_x] \partial_x(\underline{\mfT}^{-1}Q_1[\underline{U}]u_x)\\\\&&
-\varepsilon^2\mu[\Lambda^s, \underline{h}b_x^2] \underline{\mfT}^{-1}Q_1[\underline{U}]u_x 
+[\Lambda^s,Q_1[\underline{U}]]u_x.
\end{eqnarray*}
To control the term $\big(\big[\Lambda^s,Q_1[\underline{U}]\big]u_x,\Lambda^su\big)$ we use the explicit expression of $Q_1[\underline{U}]$:\\
 $$ 
 Q_1[\underline{U}]f=\frac{2}{3}\varepsilon\mu\partial(\underline{h}^3\underline{u}_xf)+\varepsilon^2\mu \underline{h}^2b_x\underline{u}_xf+\varepsilon^2\mu \underline{h}^2b_{xx}\underline{u}f,
 $$
and the fact that \\
$$
 \partial_x[\Lambda^s,f]g=[\Lambda^s, \partial_x(f \cdot)]g.
$$
Similarly to control the term $\big(\partial_x\{[\Lambda^s, \underline{h}^3]\partial_x(\underline{\mfT}^{-1}Q_1[\underline{U}]u_x)\},\Lambda^su\big)$ we use
the explicit expression of $Q_1[\underline{U}]$, the commutator estimate  (\ref{com2})  and 
Lemma \ref{proprim'}. Indeed,
\begin{eqnarray*}
\big(\partial_x\{[\Lambda^s, \underline{h}^3]\partial_x(\underline{\mfT}^{-1}Q_1[\underline{U}]u_x)\},\Lambda^su\big)&=& 
-\frac{2}{3}\varepsilon\mu\big([\Lambda^s, \underline{h}^3]\partial_x(\underline{\mfT}^{-1}\partial_x(\underline{h}^3\underline{u}_x u_x)),\Lambda^su_x\big)\\\\
&&-\varepsilon^2\mu\big([\Lambda^s, \underline{h}^3]\partial_x(\underline{\mfT}^{-1}(\underline{h}^2 b_x\underline{u}_x u_x)),\Lambda^su_x\big)\\\\
&&-\varepsilon^2\mu\big([\Lambda^s, \underline{h}^3]\partial_x(\underline{\mfT}^{-1}(\underline{h}^2 b_{xx}\underline{u}u_x)),\Lambda^su_x\big).
\end{eqnarray*}
and thus, after remarking that 
\begin{eqnarray*}
\vert \partial_x(\underline{\mfT}^{-1}\partial_x(\underline{h}^3\underline{u}_x u_x))\vert_{H^{s-1}}&\leq& \vert \underline{\mfT}^{-1}\partial_x(\underline{h}^3\underline{u}_x u_x)\vert_{H^{s}}\\ 
&\leq&\parallel \mfT^{-1}\partial_x\parallel_{H^{s}(\R)\rightarrow H^{s}(\R)}\vert \underline{h}^3\underline{u}_x u_x\vert_{H^{s}}.
\end{eqnarray*}
we can proceed as for the control of $B_3$ to get
$$
\vert B_5\vert  \leq \varepsilon C(E^s(\underline{U}))E^s(U)^2.
$$\\ 
$\bullet$ Estimate of $(\Lambda^s B(\underline{U}),S\Lambda^s U)$. Note first that 
\begin{equation*}
B(\underline{U})=\left( 
\begin{array}{c}
 \varepsilon b_x \underline{u} \\\\ 
 \underline{\mfT}^{-1} q_2(\underline{U})
\end{array}
\right)
\end{equation*}
where, $q_2(\cdot) $ as in (\ref{decompose}),
so that
\begin{eqnarray*}
(\Lambda^sB(\underline{U}),S\Lambda^s U)&=&(\Lambda^s( \varepsilon b_x \underline{u}),\Lambda^s\zeta)
+(\Lambda^s( \underline{\mfT}^{-1} q_2(\underline{U})), \underline{\mfT}\Lambda^su)\\\\
&=&(\Lambda^s( \varepsilon b_x \underline{u}),\Lambda^s\zeta)
-\big([\Lambda^s, \underline{\mfT}] \underline{\mfT}^{-1} q_2(\underline{U}) ,\Lambda^su\big)\\\\
&& +  \big(\Lambda^s q_2(\underline{U}), \Lambda^su\big). 
\end{eqnarray*}
Using again here the explicit expressions of $ \underline{\mfT}$, $q_2(\underline{U})$ and Lemma  \ref{proprim'}, we get
$$
(\Lambda^sB(\underline{U}),S\Lambda^s U) \leq \varepsilon C(E^s(\underline{U}))E^s(U).
$$
$\bullet$ Estimate of $(\Lambda^su,[\partial_t,\underline{\mfT}]\Lambda^su)$. We have that
\begin{eqnarray*}
(\Lambda^su,[\partial_t,\underline{\mfT}]\Lambda^su)&=&(\Lambda^su,\partial_t\underline{h}\Lambda^su)+\frac{\mu}{3}(\Lambda^su_x,\partial_t\underline{h}^3\Lambda^su_x)\\&&
-\frac{\varepsilon\mu}{2}(\Lambda^su,\partial_t\underline{h}^2 b_x\Lambda^su_x)-\frac{\varepsilon\mu}{2}(\Lambda^su_x,\partial_t\underline{h}^2 b_x\Lambda^su)\\&&
+\varepsilon^2\mu(\Lambda^su,\partial_t\underline{h} b_x^2 \Lambda^su).
\end{eqnarray*}
Controlling these terms by   $\varepsilon C(E^s(\underline{U}),\vert\partial_t\underline{h} \vert_{L^{\infty}})E^s(U)^2$
follows directly from a Cauchy-Schwarz inequality and an integration by 
parts.\\

Gathering the informations provided by the above estimates and using the fact that $H^s(\R)\subset W^{1,\infty}$, we get
$$
e^{\varepsilon\lambda t}\partial_t (e^{-\varepsilon\lambda t}E^s(U)^2) 
	\leq  \varepsilon\big(C(E^s(\underline{U}),\vert\partial_t\underline{h} \vert_{L^{\infty}})-\lambda\big)E^s(U)^2+\varepsilon C(E^s(\underline{U}))E^s(U).
$$
Taking $\lambda=\lambda_T$ large enough (how large depending on 
$\sup_{t\in [0,\frac{T}{\varepsilon}]}C(E^s(\underline{U}),\vert\partial_t\underline{h} \vert_{L^{\infty}})$
to have the first term of the right hand side negative for all $t\in [0,\frac{T}{\varepsilon}]$, one deduces
$$
	\forall t\in [0,\frac{T}{\varepsilon}],\qquad
	e^{\varepsilon\lambda t}\partial_t (e^{-\varepsilon\lambda t}E^s(U)^2) 
	\leq\varepsilon C(E^s(\underline{U}))E^s(U).
$$
Integrating this differential inequality yields therefore
$$
	\forall t\in [0,\frac{T}{\varepsilon}],\qquad
	E^s(U)\leq e^{\varepsilon\lambda_{T} t}E^s(U_0)+\varepsilon \int^{t}_{0} e^{\varepsilon\lambda_T( t-t')}C(E^s(\underline{U})(t'))dt'.
$$
\end{proof}
\subsection{Main result}\label{mr}
In this subsection we prove the main result of this paper, which shows
well-posedness of the Green-Naghdi equations  over large times.
\begin{theorem}\label{th1}
	Let $b\in C_b^{\infty}(\R)$,  $t_0>1/2$, $s\geq t_0+1$. Let also the initial condition
	 $U_0=(\zeta_0,u_0)^T\in X^s$,
	 and satisfy (\ref{depthcond}). Then there exists 
         a maximal $T_{max}>0$, uniformly bounded from below with respect to $\varepsilon,\mu\in (0,1)$, such that the Green-Naghdi equations (\ref{GN1}) admit 
	 a unique solution $U=(\zeta,u)^T\in X^s_{T_{max}}$ with the initial condition $(\zeta_0,u_0)^T$
         and preserving the nonvanishing depth condition (\ref{depthcond}) for any $t\in [0,\frac{T_{max}}{\varepsilon})$.
         In particular if $T_{max}<\infty$ one has
         $$ \vert U(t,\cdot)\vert_{X^s}\longrightarrow\infty\quad\hbox{as}\quad t\longrightarrow T_{max},$$
         or
         $$ \inf_\R h(t,\cdot)=\inf_\R1+\varepsilon(\zeta(t,\cdot)-b(\cdot))\longrightarrow 0 \quad\hbox{as}\quad t\longrightarrow T_{max}.$$
Morever, the following conservation of energy property holds
          $$
          \dsp\partial_t\Big(\vert\zeta\vert_2^2+(h u,u)+\mu(h\mathcal{T}u,u)\Big)=0,
          $$
where $\mathcal{T}=\mathcal{T}[h,\varepsilon b]$.
\end{theorem}
\begin{remark}
For 2D surface waves, non flat bottoms, B. A. Samaniego and D. Lannes \cite{AL2} proved  a well-posedness result to  the Green-Naghdi using a Nash-Moser scheme. Our result only use a standard Picard iterative and there is therefore no loss of regularity of the solution with respect to
the initial condition. In the one-dimensional case and for flat bottoms, our result coincides
with the one proved by Li in \cite{li}.
\end{remark}
\begin{remark}
Our approach does not admit a straightforward generalization to the 2D case. The main 
reason is that
the natural energy norm $X^s$ is then given by
$$
\vert U\vert_{X^s}^2=\vert\zeta\vert_{H^s}^2+\vert u\vert_{H^s}^2
+\mu \vert \nabla\cdot u\vert_{H^s}^2,
$$
which does \emph{not} control the $H^1(\R^2)$ norm of $u$ (since $u$ takes its values in $\R^2$, the information on the rotational of $u$ is missing).
\end{remark}
\begin{remark}
No smallness assumption on $\varepsilon$ nor $\mu$ is required in the theorem. The fact that $T_{max}$
is uniformly bounded from below with respect to these parameters allows us to say that if some
smallness assumption is made on $\varepsilon$, then the existence time becomes larger, namely of order $O(1/\varepsilon)$. This  is consistent with the existence time obtained for the (simpler)  physical models derived under some smallness assumption on $\varepsilon$, like the Boussinesq models. In fact, such models can be derived from the Green-Naghdi equations \cite{LB}. The present theorem also
has some direct implication for the justification of variable-bottom Camassa-Holm equations \cite{Israwi}.
\end{remark}
\begin{proof}
We want to construct a sequence of approximate solution $(U^n=(\zeta^n,u^n))_{n\ge 0}$ by the 
induction relation
\begin{equation}\label{approximatesys}
         U^0=U_0,\quad\mbox{ and }\quad
	\forall n\in\N, \quad
	\left\lbrace
	\begin{array}{l}
	\dsp\partial_t U^{n+1}+A[U^n]\partial_x U^{n+1}+B (U^n)=0;
        \\
	\dsp U^{n+1}_{\vert_{t=0}}=U_0.
	\end{array}\right.
\end{equation}
By Proposition \ref{ESprop}, we know that there is a unique solution $U^{n+1}\in C([0,\infty);X^s)$ to
(\ref{approximatesys}) if $U^{n}\in C([0,\infty);X^s)$ and $U^n$ satisfies (\ref{depthcond}) for all times.
Let $R>0$ be such that $E^s(U^0)\leq R/2$,
it follows from Proposition \ref{ESprop} that $U^{n+1}$ satisfies the following inequality 
\begin{eqnarray*}
	E^s(U^{n+1}(t))&\leq&e^{\varepsilon\lambda_T t}E^s(U^0)+ \varepsilon \int^{t}_{0} e^{\varepsilon\lambda_T( t-t')}C(E^s(U^n(t'))dt',
\end{eqnarray*}
we suppose now that $$\sup_{t\in [0,\frac{T}{\varepsilon}]}E^s(U^n(t))\leq R,$$ therefore
$$
E^s(U^{n+1}(t))\leq R/2+(e^{\varepsilon\lambda_Tt}-1)(R/2+\frac{C(R) }{\lambda_T}). 
$$
Hence, there is $T>0$ small enough such that
$$
\sup_{t\in [0,\frac{T}{\varepsilon}]}E^s(U^{n+1}(t))\leq R.
$$
Using now the link between $E^s(U)$ and $\vert U\vert_{X^s}$ given by Lemma \ref{lemmaes} 
we get 
$$
\sup_{t\in [0,\frac{T}{\varepsilon}]}\vert U^{n+1}(t)\vert_{X^s}\leq C\big(\frac{1}{h_0}\big)R.
$$
We also know from the equations  that
\begin{equation*}
\partial_t \zeta^{n+1}=-h^{n}u^{n+1}_x-\varepsilon\zeta^{n}_xu^{n+1}+\varepsilon b_x u^{n+1}.
\end{equation*}
Hence, one gets 
 \begin{equation}\label{zeta_t}
\vert \partial_t h^{n+1}\vert_{L^{\infty}}=\varepsilon \vert \partial_t \zeta^{n+1}\vert_{L^{\infty}}\leq \varepsilon C\big(\frac{1}{h_0}\big) R.
\end{equation}
Since moreover
$$
h^{n+1}=h^{n+1}_{t=0}+\int_{0}^{t}\partial_t\zeta^{n+1},
$$
we can deduce from (\ref{zeta_t}) and the fact that $h^{n+1}_{t=0}=1+\varepsilon(\zeta_0-b)\ge h_0$ that it is possible to choose  $T$ small enough
for $U^{n+1}$ to satisfy  (\ref{depthcond}) on $[0,\frac{T}{\varepsilon}]$,
with $h_0$ replaced by $h_0/2$.  \\
Finally, we deduce that the Cauchy problem
\begin{equation*}
	\left\lbrace
	\begin{array}{l}
	\dsp\partial_t U^{n+1}+A[U^n]\partial_x U^{n+1}+B (U^n)=0;
        \\
	\dsp U^{n+1}_{\vert_{t=0}}=U_0
	\end{array}\right.
\end{equation*}
has a unique solution $U^{n+1}$ satisfing (\ref{depthcond}) and the inequality
\begin{eqnarray*}
	E^s(U^{n+1})&\leq&e^{\varepsilon\lambda_Tt}E^s(U_0)+ \varepsilon \int^{t}_{0} e^{\varepsilon\lambda_T( t-t')}C(E^s(U^n)(t'))dt',
\end{eqnarray*}
when $0\leq t\leq\frac{T}{\varepsilon}$ and
 $\lambda_T$ depending only  on $\sup_{t\in [0,\frac{T}{\varepsilon}]} E^s(U^n)$.
 Thanks to this energy estimate, one can conclude classically
(see e.g. \cite{AG})
to the existence of 
$$
	T_{max}=T(E^s(U_0))>0,
$$ 
and of
a unique solution $U\in X^s_{T_{max}}$ to
(\ref{GN1}) preserving the inequality (\ref{depthcond}) for any $t\in [0,\frac{T_{max}}{\varepsilon}]$  as a limit of the iterative scheme
$$
	U^0=U_0,\quad\mbox{ and }\quad
	\forall n\in\N, \quad
	\left\lbrace
	\begin{array}{l}
	\dsp\partial_t U^{n+1}+A[U^n]\partial_x U^{n+1}+B (U^n)=0;
        \\
	\dsp U^{n+1}_{\vert_{t=0}}=U_0.
	\end{array}\right.
$$
The fact that $T_{max}$ is bounded from below by some $T>0$ independent of $\varepsilon,\mu\in (0,1)$ follows from the analysis above, 
while the behavior of the solution as $t\to T_{max}$ if $T_{max}<\infty$ follows from standard
continuation arguments.\\
Though the conservation of the energy can be found in some references (e.g. \cite{Chazel}), we reproduce it
in Appendix \ref{sectconserv} for the sake of completeness.

\end{proof}

\appendix

\section{Existence of solutions for the linearized equations}\label{appendix}
In this section we examine existence, uniqueness, and regularity for solutions to the following system of equations:
\begin{equation}\label{GNlsysapp}
	\left\lbrace
	\begin{array}{l}
	\dsp\partial_t U+A[\underline{U}]\partial_x U=f;
        \\
	\dsp U_{\vert_{t=0}}=U_0,
	\end{array}\right.
\end{equation}
where $\underline{U}=(\underline{\zeta}, \underline{u})^T$ $\in X^{s}_{T}
$ is such that $\partial_t \underline{U} \in X^{s-1}_{T}$ 
and satisfy the condition  (\ref{depthcond}) on $[0,\frac{T}{\varepsilon}]$. 
We begin the proof by the following lemma (see for instance \cite{taylor}): 
\begin{lemma}\label{prop phi}
Let $\varphi\in C^{\infty}_{0}(\R)$, such that $\varphi(r)=1$ for $\vert r\vert\leq 1$. Let also
$$ J^{\delta}=\varphi(\delta\vert D\vert), \quad \delta >0.$$
Then:\\
(i)  $\forall s, s' \in \R$, $ J^{\delta}$: $H^{s}(\R)\longmapsto  H^{s'}(\R)$ is a bounded linear operator.\\
(ii) $J^{\delta}$ commutes with $\Lambda^s$ and is self-adjoint operator.\\ 
(iii) $\forall f \in C^1(\R)\cap L^\infty(\R)$, $v\in L^2(\R)$ there exists $C$ independent of $\delta$ such that
$$
\vert[f,J^{\delta}]v\vert_{H^1}\leq C \vert f\vert_{C^1}\vert v\vert_{L^2}.
$$
(iv) $\forall f \in H^s(\R)$, $s>\frac{1}{2}$, 
 $ J^{\delta}f \in L^{\infty}(\R)$ with 

$$\vert J^{\delta} f\vert _{L^{\infty}}\leq C \vert f\vert_{L^{\infty}}$$ 
where $C$ is a constant independent of $\delta$.
\end{lemma}
Our strategy will be to obtain a solution to (\ref{GNlsysapp}) as a limit of solutions $U_{\delta}$ to
\begin{equation}\label{GNlsysapproch}
	\left\lbrace
	\begin{array}{l}
	\dsp\partial_t U_{\delta}+J^{\delta}A[\underline{U}] J^{\delta}\partial_xU_{\delta}=f;
        \\
	\dsp U_{\delta_{\vert_{t=0}}}=U_0.
	\end{array}\right.
\end{equation}
For any $\delta >0$, $A^{\delta}=J^{\delta}A[\underline{U}] J^{\delta}$ is  a bounded linear operator on each 
$X^s$, and $F^{\delta}=A^{\delta}\partial_x+B (\underline{U})\in C^1(X^s)$ so by 
Cauchy-Lipschitz the ODE (\ref{GNlsysapproch})  has a unique solution, $U_{\delta}\in C([0,T/\varepsilon],X^s)$. 
Our task will be to obtain estimates on $U_{\delta}$, independent 
of $\delta\in (0,1)$ and to show that the solution $U_{\delta}$ has a limit as $\delta\searrow 0$ solving (\ref{GNlsysapp}).
To do this, we remark that
\begin{eqnarray}\label{qtsctrapp}
\frac{1}{2}\partial_t E^s(U_{\delta})^2 &=&
-(SA[\underline{U}]\Lambda^s\partial_x J^\delta U_\delta,\Lambda^s J^\delta U_\delta)- \big(\big[\Lambda^s,A[\underline{U}]\big]\partial_x J^\delta U_\delta,S\Lambda^s J^\delta U_\delta\big)\nonumber\\\nonumber\\&&
+(\Lambda^s f,S\Lambda^s U_\delta)+\frac{1}{2}(\Lambda^su_\delta,[\partial_t,\underline{\mfT}]\Lambda^su_\delta)\nonumber\\\nonumber\\&&
+\big(\big[S,J^\delta\big]\Lambda^s U_\delta,A[\underline{U}]\Lambda^s \partial_x J^\delta U_\delta\big)+\big(\big[S,J^\delta\big]\Lambda^s U_\delta, \big[\Lambda^s,A[\underline{U}]\big]\partial_x J^\delta U_\delta\big).
\end{eqnarray}
Note that we do not give any details for the control of the components of  the r.h.s (\ref{qtsctrapp}) other than the last two terms  because the others  can  be handled exactly as in Proposition \ref{ESprop}.
To estimate the last two terms of the r.h.s  (\ref{qtsctrapp}), we have that
\begin{eqnarray*}
 \big(\big[S,J^\delta\big]\Lambda^s U_\delta,A[\underline{U}]\Lambda^s \partial_x J^\delta U_\delta\big)&=&
\big(\big[\underline{\mfT},J^\delta\big]\Lambda^s u_\delta,\underline{\mfT}^{-1}(\underline{h}\Lambda^s \partial_x J^\delta u_\delta)\big)\\\\
&&+\big(\big[\underline{\mfT},J^\delta\big]\Lambda^s u_\delta,\varepsilon\underline{u}\Lambda^s \partial_x J^\delta u_\delta\big)\\\\&&
+\big(\big[\underline{\mfT},J^\delta\big]\Lambda^s u_\delta,\underline{\mfT}^{-1}Q_1[\underline{U}]\Lambda^s \partial_x J^\delta u_\delta\big).
\end{eqnarray*}
One can check by using the explicit expression of $\underline{\mfT}$ that
\begin{eqnarray*}
 \big[\underline{\mfT},J^\delta\big]&=& \big[\underline{h},J^\delta\big]-\frac{\mu}{3}\partial_x\big[\underline{h^3},J^\delta\big]\partial_x
-\frac{\varepsilon\mu}{2}\big[\underline{h^2}b_x,J^\delta\big]\partial_x\\\\&&
+\frac{\varepsilon\mu}{2}\partial_x\big[\underline{h}^2b_x,J^\delta\big]
+\varepsilon^2\mu\big[\underline{h}b_x^2,J^\delta\big].
\end{eqnarray*}
One deduces directly from Lemma \ref{prop phi}, an integration by parts, the
 Cauchy-Schwarz inequality and the explicit expression of $Q_1[\underline{U}]$ that
$$
\big(\big[S,J^\delta\big]\Lambda^s U_\delta,A[\underline{U}]\Lambda^s \partial_x J^\delta U_\delta\big)
\leq C E^s( U_\delta)^2,
$$
similarly, one can conclude
$$
\big(\big[S,J^\delta\big]\Lambda^s U_\delta, \big[\Lambda^s,A[\underline{U}]\big]\partial_x J^\delta U_\delta\big)
\leq C E^s( U_\delta)^2,
$$
where $C$ is a constant independent of $\delta$. By the Proposition \ref{ESprop}, we have 
\begin{eqnarray*}
 &&-(SA[\underline{U}]\Lambda^s\partial_x J^\delta U_\delta,\Lambda^s J^\delta U_\delta)- \big(\big[\Lambda^s,A[\underline{U}]\big]\partial_x J^\delta U_\delta,S\Lambda^s J^\delta U_\delta\big)\nonumber\\\\&&
+(\Lambda^s f,S\Lambda^s U_\delta)+\frac{1}{2}(\Lambda^su_\delta,[\partial_t,\underline{\mfT}]\Lambda^su_\delta)\\\\
&&\leq C E^s( U_\delta)^2+C E^s(f)^2.
\end{eqnarray*}
Consequently, we obtain an estimate of the form
\begin{equation}\label{estiamtG}
\frac{d}{dt}E^s( U_\delta)^2\leq C E^s( U_\delta)^2+C E^s(f)^2.
\end{equation}
Thus Gronwall's inequality yields an estimate 
\begin{equation}\label{estiamtG}
E^s( U_\delta)^2\leq C(t)\big[E^s( U_0)^2+\sup_{[0,t]}E^s(f)^2\big],
\end{equation}
independent of $\delta\in (0,1)$.
Thanks to this energy estimate, one can conclude classically
(see e.g. \cite{taylor})
to the existence of
a unique solution $U\in C([0,T], X^s)$ to (\ref{GNlsysapp}).

\section{Conservation of the energy}\label{sectconserv}

In order to prove that
  $$
          \dsp\partial_t\Big(\vert\zeta\vert_2^2+(h u,u)+\mu(h\mathcal{T}u,u)\Big)=0,
  $$
we multiply the first equation of (\ref{GN1}) by $\zeta$ and the second by $u$, 
integrateon $\R$, and sum both equations to find
$$
\frac{1}{2}\partial_t\vert\zeta\vert_2^2+(\partial_x(hu),\zeta)+(\partial_t u,hu)+\mu(h\mathcal{T}\partial_t u,u)+(\partial_x\zeta,hu)
+\varepsilon(u\partial_xu,hu)+\mu\varepsilon(\mathcal{Q} u,u)=0.
$$
Therefore
$$
\frac{1}{2}\partial_t\vert\zeta\vert_2^2+\frac{1}{2}\partial_t(hu,u) - \frac{1}{2}( \partial_th,u^2)+\varepsilon(u\partial_xu,hu)+\mu(h\mathcal{T}\partial_t u,u)
+\mu\varepsilon(\mathcal{Q} u,u)=0,
$$
where the term $\mathcal{Q}u$ is defined as:
\begin{eqnarray*}
 \mathcal{Q} u&=& -\frac{1}{3}\partial_x[(h^3(u\partial_x^2 u-(\partial_x u)^2)+\frac{\varepsilon}{2}[\partial_x(h^2u\partial_x(u\partial_xb)-h^2\partial_xb(u\partial_x^2 u-(\partial_x u)^2)]\\&&
+\varepsilon^2h\partial_xb(u\partial_x(u\partial_xb)).
\end{eqnarray*}
Using now  the fact that $h=1+\varepsilon(\zeta-b)$ and the first equation of (\ref{GN1}), we get
\begin{eqnarray*}
-\frac{1}{2}(\partial_th,u^2)+\varepsilon(u\partial_xu,hu)=-\frac{\varepsilon}{2}(\partial_x(hu),u^2)+\varepsilon(u\partial_xu,hu)=0.
\end{eqnarray*}
Thus, 
\begin{equation}\label{CE}
\frac{1}{2}\partial_t\vert\zeta\vert_2^2+\frac{1}{2}\partial_t (hu,u)+\mu(h\mathcal{T}\partial_t u,u)+\mu\varepsilon(\mathcal{Q} u,u)=0.
\end{equation}
Regarding now the term $\mu(h\mathcal{T}\partial_t u,u)$, we remark as in \cite{Chazel} that
\begin{eqnarray*}
\mu(h\mathcal{T}\partial_t u,u)&=&\mu(\mathcal{T}_1^*h\mathcal{T}_1\partial_t u,u)+\mu(\mathcal{T}_2^*h\mathcal{T}_2\partial_t u,u),\\
&=&\mu(h\mathcal{T}_1\partial_t u,\mathcal{T}_1u)+\mu(h\mathcal{T}_2\partial_t u,\mathcal{T}_2u),\\
&=&\mu(h(\partial_t(\mathcal{T}_1 u)-\partial_t\mathcal{T}_1 u),\mathcal{T}_1u)+\mu(h\partial_t (\mathcal{T}_2u),\mathcal{T}_2u),
\end{eqnarray*}
with $\mathcal{T}_j^*(j=1,2)$ denoting the adjoint of the operators 
$\mathcal{T}_j$ given by
$$
\mathcal{T}_1u=\frac{h}{\sqrt{3}}\partial_xu-\varepsilon\frac{\sqrt{3}}{2}\partial_xbu, \quad \hbox{and} \quad \mathcal{T}_2u=\frac{\varepsilon}{2}\partial_xbu.
$$ 
It comes:
\begin{eqnarray*}
\mu(h\mathcal{T}\partial_t u,u)&=&\frac{\mu}{2}\partial_t(h\mathcal{T}_1 u,\mathcal{T}_1u)-\frac{\mu}{2}(\partial_th,(\mathcal{T}_1u)^2)
-\mu(h(\partial_t\mathcal{T}_1 u,\mathcal{T}_1u)\\&&+\frac{\mu}{2}\partial_t(h\mathcal{T}_2 u,\mathcal{T}_2u)-\frac{\mu}{2}(\partial_th,(\mathcal{T}_2u)^2),\\
&=&\frac{\mu}{2}\partial_t(h\mathcal{T}_1 u,\mathcal{T}_1u)+\frac{\mu}{2}\partial_t(h\mathcal{T}_2 u,\mathcal{T}_2u)-\frac{\mu}{2}(\partial_th,(\mathcal{T}_1u)^2+
(\mathcal{T}_2u)^2)\\&&-\mu(h\partial_t\mathcal{T}_1 u,\mathcal{T}_1u).
\end{eqnarray*}
Inject this result in (\ref{CE}) to get:
\begin{eqnarray*}
 \dsp\frac{1}{2}\partial_t\Big(\vert\zeta\vert_2^2+(h u,u)+\mu(h\mathcal{T}u,u)\Big)&=&\frac{\mu}{2}(\partial_th,(\mathcal{T}_1u)^2+
(\mathcal{T}_2u)^2)\\&&+\mu(h\partial_t\mathcal{T}_1 u,\mathcal{T}_1u)-\mu\varepsilon(\mathcal{Q} u,u)\nonumber.
\end{eqnarray*}
Noting that $h\partial_t\mathcal{T}_1 u=\partial_th(\mathcal{T}_1 u+\sqrt{3}\mathcal{T}_2 u)$, it comes:
\begin{eqnarray*}
\frac{1}{2}(\partial_th,(\mathcal{T}_1u)^2+
(\mathcal{T}_2u)^2)+(h\partial_t\mathcal{T}_1 u,\mathcal{T}_1u)&=&\frac{1}{2}\Big(\partial_th,3(\mathcal{T}_1u)^2+
(\mathcal{T}_2u)^2+2\sqrt{3}\mathcal{T}_1u\mathcal{T}_2u\Big),\\
&=&\frac{\mu}{2}\Big(\partial_th,(\sqrt{3}\mathcal{T}_1u+
\mathcal{T}_2u)^2\Big),\\
&=&\varepsilon\Big(u,\frac{h}{2}\partial_x(\sqrt{3}\mathcal{T}_1u+
\mathcal{T}_2u)^2\Big),
\end{eqnarray*}
where we used here the first equation of (\ref{GN1}). Finally, we get:
\begin{eqnarray*}
 \dsp\frac{1}{2}\partial_t\Big(\vert\zeta\vert_2^2+(h u,u)+\mu(h\mathcal{T}u,u)\Big)&=&\mu\varepsilon\Big(\frac{h}{2}\partial_x(\sqrt{3}\mathcal{T}_1u+
\mathcal{T}_2u)^2-\mathcal{Q} u,u\Big)\nonumber.
\end{eqnarray*}
One can easily show that $\Big(\dsp\frac{h}{2}\partial_x(\sqrt{3}\mathcal{T}_1u+
\mathcal{T}_2u)^2-\mathcal{Q} u,u\Big)=0$,  which implies easily the result.

\subsection*{Acknowledgments}\ The author is grateful to David Lannes for encouragement and many helpful discussions.

\providecommand{\href}[2]{#2}
\end{document}